\input gtmacros.tex
\input amsnames.tex
%
%
\medskipamount=6pt plus 2pt minus 2pt
\def\cS{{\cal S}}
\def\cee{{\Bbb C}}
\def\real{{\Bbb R}}
\def\zed{{\Bbb Z}}
\def\Bspin{B\hbox{\rm spin}}
\def\endrk{\medskip}
\def\cite#1{{\rm [#1]}}
\def\frac#1#2{{#1\over#2}}
\mathsurround0.7pt
\count0=41
\reflist 

\refkey\DK
{\bf S Donaldson}, {\bf P Kronheimer},
{\it The geometry of four-manifolds},
Oxford Mathematical Monographs, 
Clarendon Press, 
Oxford, UK (1990)

\refkey\HH
{\bf F Hirzebruch}, {\bf H Hopf}, 
{\it Felder von Fl\"achenelementen in 4--dimension\-alen 
Mannigfaltigkeiten}, 
Math. Ann. 136 (1958) 156--172

\refkey\KS
{\bf R Kirby}, {\bf L Siebenmann},
{\it Foundational essays on topological manifolds,\break smoothings and
triangulations}, 
Ann. Math. Studies  88, 
Princeton Univ. Press,
Princeton, NJ (1977) 

\refkey\Mi
{\bf J Milnor}, 
{\it On spaces having the homotopy type of a $CW$--complex}, 
Trans. Amer. Math. Soc.  90 (1959) 272--280

\refkey\Mil
{\bf J Milnor}, 
{\it Spin structures on manifolds}, 
L'Ensignement Math. 8 (1962)  198--203

\refkey\Mo
{\bf J Morgan}, 
{\it The Seiberg--Witten equations and applications to the topology
of smooth 4--manifolds}, 
Princeton Univ. Press, 
Princeton, NJ (1996) 

\refkey\S 
{\bf M Spivak},  
{\it Spaces satisfying Poincar\'e duality},  
Topology  6 (1967)  77--102 

\refkey\W
{\bf G Whitehead}, 
{\it Elements of homotopy theory}, 
Graduate Texts in Math.  61, 
Springer-Verlag, 
Berlin (1978) 

\endreflist 
\input gtoutput
\volumenumber{1}\papernumber{4}\volumeyear{1997}
\pagenumbers{41}{50}
\published{23 October 1997}
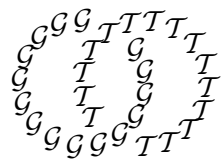
\title{Spin$^c$--structures and homotopy equivalences}

\author{Robert E Gompf} 
\address{Department of Mathematics\\The University of Texas at 
Austin\\Austin, TX 78712-1082 USA}
\email{gompf@math.utexas.edu}

\abstract
We show that a homotopy equivalence between manifolds induces
a correspondence between their spin$^c$--structures, even in the presence
of 2--torsion.  This is proved by generalizing spin$^c$--structures to
Poincar\'e complexes.  A procedure is given for explicitly computing
the correspondence under reasonable hypotheses.
\endabstract

\proposed{Ronald Stern}\received{16 May  1997}
\seconded{Robion Kirby, Dieter Kotschick}\accepted{17 October 1997}

\maketitlepage

%
\def\ThmFive{5}
\def\PropOne{1}

\section{Introduction}

The theory of spin$^c$--structures has attained new importance through its 
recent application to the topology of smooth 4--manifolds. 
Among smooth, closed, oriented 4--manifolds (with $b_1+b_+$ odd) a typical 
homeomorphism type contains many diffeomorphism types. 
The only invariants known to distinguish such diffeomorphism types are 
those arising from gauge theory, as pioneered by Donaldson (eg \cite{\DK}). 
The most efficient approach currently known is to assign a 
{\it Seiberg--Witten\/} invariant (eg \cite{\Mo}) to any such 4--manifold 
$X$ with a fixed spin$^c$--structure. 
To extract the most information from
these invariants, one must understand how spin$^c$--structures 
transform under homeomorphisms. 
This is straightforward if $H^2(X;\zed)$ has no 2--torsion (for example, 
if $X$ is simply connected), for then the Chern class will distinguish any 
two spin$^c$--structures on $X$. 
The general case is less obvious, however. 
In high dimensions, a homeomorphism between smooth manifolds need not be 
covered by an isomorphism of their tangent bundles. 
While such isomorphisms always exist in dimension~4, they are not canonical, 
and automorphisms of the tangent bundle covering $id_X$ may permute the 
spin$^c$--structures on $X$. 
(For example, such an automorphism over $\real P^3$ or $\real P^3\times S^1$ 
can be constructed from the diffeomorphism $\real P^3 \to SO(3)$.)  
In this note, we show how to canonically assign to any orientation-preserving 
proper homotopy equivalence $X_1\to X_2$ between manifolds
a correspondence between 
spin$^c$--structures on $X_1$ and those on $X_2$. 

Our approach is to generalize the theory of spin and spin$^c$--structures 
from $SO(n)$ to more general structure groups $H$. 
Most of the homotopy  of $SO(n)$ does not enter into the theory. 
In fact, it suffices for $H$ to be path connected with a nontrivial double 
cover so that we can generalize the definition $\hbox{spin}^c(n) = 
(\hbox{spin}(n)\times \hbox{spin}(2))/\zed_2$. 
The resulting theory generalizes the classical theory in the obvious way, 
for example, with spin$^c$--structures on a bundle $\xi$ over $X$
classified by $H^2(X;\zed)$ whenever 
$W_3 (\xi)=0$ (Proposition~\PropOne). 
Ultimately, the map $BSO \to BSG$ of classifying spaces allows us to 
generalize spin$^c$--structures from smooth manifolds to Poincar\'e complexes, 
and the latter theory has the required functoriality with respect to homotopy 
equivalences by naturality of the Spivak normal fibration 
(Theorem \ThmFive).		
Under reasonable hypotheses, one can explicitly compute the correspondence 
of spin$^c$--structures induced by a homotopy equivalence; a procedure is 
given following Theorem~\ThmFive. 
The concluding remarks include other characterizations of classical 
spin$^c$--structures. 

\rk{Acknowledgements}The author wishes to thank Larry Taylor for 
helpful discussions, and 
the Mathematical Sciences Research Institute for their hospitality. 
The author is partially supported by NSF grant DMS-9625654. Research at MSRI
is supported in part by NSF grant DMS-9022140.

\section{Generalized spin$^c$--structures}

A naive approach to generalizing the theory of spin and spin$^c$--structures
would be to define $\hbox{spin}(H)$ to be a preassigned double cover of a
path connected topological group $H$, and let $\hbox{spin}^c (H)$ denote the
group $\hbox{spin}(H \times SO(2))$ diagonally double covering $H\times SO(2)$.
One could then generalize the theory in the obvious way, using principal
$\hbox{spin}(H)$ and $\hbox{spin}^c (H)$--bundles, the natural epimorphisms
from $\hbox{spin}^c (H)$ to $H$ and $SO(2)$, and the involution of 
$\hbox{spin}^c(H)$ induced by conjugation on $SO(2) = U(1)$.  However, to
avoid the difficulties of adapting principal bundle theory to spherical
fibrations, we translate the argument into the language of classifying
spaces, replacing epimorphisms of groups with kernel $\zed_2$ or $ SO(2)$ by
fibrations of the corresponding
classifying spaces with fiber $B\zed_2 = K(\zed_2,1) = \real P^\infty$
or $BSO(2) = K(\zed,2) = \cee P^\infty$, respectively.  We remove the groups
from the theory while keeping the suggestive notation, obtaining a
theory of spin and spin$^c$--structures on bundles or 
fibrations classified by a universal bundle (fibration) $\xi_H\to BH$, 
where $BH$ is homotopy equivalent to a simply connected $CW$--complex, and 
a nonzero class $w\in H^2 (BH;\zed_2)$ is specified (corresponding to a choice
of double cover of $H$). We can recover the classical theory by
setting $BH=BSO(n)$ ($n\ge2$), with $w$ the unique nonzero class $w^*\in H^2
(BSO(n);\zed_2)\cong\zed_2$. 

Recall \cite{\W} that any map $f\co X\to Y$ can be transformed into a fibration
by replacing $X$ by the space $P$ of paths from $X$ to $Y$ in the mapping
cylinder of $f$.  The initial point fibration $p_0\co P\to X$ has contractible
fiber, and the endpoint fibration $p_1\co P\to Y$ is homotopic to $f\circ p_0$.
The fiber $F$ of $p_1$ is homotopy equivalent to a $CW$--complex if $X$ 
and $Y$ are \cite{\Mi}, and $p_0|F$ is a fibration with fiber the loop
space $\Omega Y$.

Now let $ (BH,w)$ be as above.  Then
$w$ defines epimorphisms $H_2(BH;\zed_2)\to\zed_2$ 
and hence $\varphi_w\co \pi_2 (BH)\to \zed_2$. We apply the previous paragraph
to the map $BH\to K(\zed_2,2)$ induced by $\varphi_w$, and let 
$\Bspin(H,w)$ denote the fiber $F$. 
The fibration $\Bspin (H,w)\to BH$ induces isomorphisms of 
$\pi_i (\Bspin (H,w))$ with $\ker \varphi_w$ for $i=2$ and $\pi_i(BH)$ 
otherwise, and its fiber is 
$K(\zed_2,1)= \real P^\infty$. 
Now we define $\Bspin^c(H,w)$ to be $\Bspin (H^\times ,w+w^*)$, where 
$BH^\times = BH\times BSO(2)$. 
We immediately obtain fibrations $p_H$ and $p_{SO(2)}$ of 
$\Bspin^c (H,w)$ over $BH$ and $BSO(2)$, whose fibers are
$\Bspin (SO(2),w^*) = K(\zed,2)$ and 
$\Bspin (H,w)$, respectively, and each fibration restricted to 
the opposite fiber is the map arising from the
definition of $\Bspin( \cdot )$. (Compare with the projections of
$\hbox{spin}^c(H,w)$ to $H$ and $SO(2)$ on the level of groups.)  By 
obstruction theory, complex conjugation on 
the second factor $BSO(2)=\cee P^\infty$ of $BH^\times$ 
lifts uniquely from $BH^\times$ to a map on $\Bspin^c(H,w)$  
whose square is fiber homotopic to the identity, and the map
is homotopic to conjugation on each 
$\cee P^\infty$--fiber of $p_H$. 

To define spin$^c$--structures over $H$, recall that an $H$--bundle (or fibration)
$\xi\to X$ over a $CW$--complex is classified by a bundle map 
\def\strt{\vrule height 0pt width 0pt depth 5pt}
$$\matrix{ 
\xi\strt&\buildrel {\tilde f}\over\longrightarrow & \xi_H\cr 
\Big\downarrow&&\Big\downarrow\cr
X&\buildrel f\over\longrightarrow&BH.\cr} 
$$
For two choices of classifying map $\tilde f$, there is a canonical
homotopy (up to homotopy rel 0,1) between the corresponding maps $f$,
characterized by lifting to a homotopy of the maps $\tilde f$
through bundle maps.  This allows us to define
spin$^c$--structures in a manner independent of the choice of $\tilde f$.

\rk{Definition}%
A {\it spin structure\/} on an $H$--bundle (fibration)
$\xi$ (relative to $w$) is a 
function assigning to each classifying bundle map $\tilde f\co \xi \to \xi_H$
a homotopy class of lifts $\hat f\co X\to\Bspin (H,w)$ of $f\co X\to BH$, 
such that for two choices 
of $\tilde f$ the canonical homotopy between the maps $f$ 
lifts to a homotopy of the corresponding maps $\hat f$. 
A {\it spin$^c$--structure\/} is defined similarly with spin replaced by 
spin$^c$.
\endrk

We denote the sets of spin and spin$^c$--structures on an $H$--bundle $\xi$ 
by $\cS(\xi,w)$ and $\cS^c(\xi,w)$, respectively. 
Note that in either case, any lift of a single $f$ 
with a specified $\tilde f$ uniquely determines such a structure,
but changing $\tilde f$ with $f$ fixed may 
result in an automorphism of $\cS(\xi,w)$ or $\cS^c(\xi,w)$. 

To define characteristic classes, let $Y\subset X$ be a possibly empty 
subcomplex,\break and let $\tau$ be a trivialization of $\xi|Y$. 
Then we can assume that the classifying map $f\co X\to BH$ of $\xi$ is 
constant on $Y$, and that $\tau$ determines the restriction $\tilde f|Y\co 
\xi|Y\to \xi_H$. 
Set $w_2 (\xi,\tau) = f^* (w) \in H^2(X,Y;\zed_2)$ and 
$W_3 (\xi,\tau) = \beta w_2 (\xi,\tau) \in H^3 (X,Y;\zed)$, where $\beta$ 
is the Bockstein homomorphism. 
Any spin$^c$--structure $s\in \cS^c(\xi,w)$ determines a homotopy class of lifts 
$\hat f\co X\to \Bspin^c (H,w)$ of $f$, and we define a trivialization $\hat\tau$ 
of $s|Y$ over $\tau$ to be a choice of $\hat f$ (within the given homotopy
class) that is constant on Y, up to homotopies through such maps.
(Equivalently, $\hat\tau$ is a spin$^c$--structure on $X/Y$ that pulls back to
$s$ on $X$.)
We define Chern classes by setting $c_1(s,\hat\tau) = \hat f^* p_{SO(2)}^* 
(c) \in H^2(X,Y;\zed)$, where $c\in H^2(BSO(2),\zed) \cong \zed$ is the 
generator $c_1(\xi_{SO(2)})$.  
If $Y$ is empty, we use the notation $w_2(\xi)$, $W_3(\xi)$, $c_1(s)$.

\proclaim{Proposition 1}		
{The set $\cS(\xi,w)$ of spin structures on an $H$--bundle (or fibration)
$\xi\to X$ is 
nonempty  if and only if $w_2(\xi)=0$. 
If so, then $H^1(X;\zed_2)$ acts freely and transitively on $\cS(\xi,w)$. 
The set $\cS^c(\xi,w)$ is nonempty if and only if $W_3(\xi)=0$, and if so, 
then $H^2(X;\zed)$ acts freely and transitively on it. 
For $s\in \cS^c(\xi,w)$ and $a\in H^2(X;\zed)$, we have $c_1(s+a)= c_1(s)+2a$. 
Conjugation induces an involution on $\cS^c(\xi,w)$ that reverses signs 
of Chern classes and the $H^2(X;\zed)$--action. 
For $Y\subset X$ and $\hat\tau$ as above, $c_1(s,\hat\tau)$ reduces 
modulo~2 to $w_2(\xi,\tau)$.}
\endproc

Thus, choosing a base point in $\cS (\xi,w)$ or $\cS^c (\xi,w)$ (if  nonempty) 
identifies it with $H^1(X;\zed_2)$ or $H^2(X;\zed)$. 

\prf
The first two sentences are immediate from obstruction theory, since the 
fiber of $\Bspin (H,w)\to BH$ is $K(\zed_2,1)$. 
In fact, $w_2(\xi,\tau)$ is the obstruction to lifting $f$ to a map 
$\hat f\co X\to\Bspin (H,w)$ with $\hat f|Y$ constant, as can be seen  by 
first considering the case where $Y$ contains the 1--skeleton of $X$. 
Similarly, $H^2(X;\zed)$ acts as required on $\cS^c(\xi,w)$ (when nonempty) 
via difference classes, since the fiber of $p_H$ is $K(\zed,2)$. 
Now recall that $\Bspin^c(H,w) = \Bspin (H^\times,w+ w^*)$ with 
$BH^\times = BH\times BSO(2)$. 
Thus, a lift of $f$ to $\hat f\co X\to\Bspin^c (H,w)$ with $\hat f|Y$ constant 
is the same as a choice of complex line bundle $L\to X$ with a 
trivialization $\tau_L$ over $Y$, together with a spin structure on 
the bundle $\xi\times L\to X$ (classified by $BH\times BSO(2)$) whose 
defining lift $\hat f$ is constant on $Y$. 
The resulting spin$^c$--structure $s$ with trivialization $\hat\tau$ over 
$\tau$  will satisfy $c_1(s,\hat \tau) = c_1(L,\tau_L)$, since $p_{SO(2)}  
\circ \hat f$ is the classifying map of $L$. 
Such a structure exists if and only if $0=w_2(\xi\times L,\tau\times \tau_L) 
= w_2(\xi,\tau) + w_2(L,\tau_L)$, or equivalently $w_2(\xi,\tau) = w_2 (L,
\tau_L) = c_1(L,\tau_L)|_2$. 
Thus, $\cS^c(\xi,w)$ is nonempty if and only if $w_2(\xi)$ has a lift to  
$H^2(X;\zed)$, ie  $W_3(\xi)=0$, and any $c_1(s,\hat\tau)$ reduces 
mod~2 to $w_2(\xi,\tau)$. 
Given $s,s'\in \cS^c(\xi,w)$, the difference class $d(s,s')$ takes 
coefficients in $\pi_2(K(\zed,2))$, where $K(\zed,2)$ is the fiber of $p_H$.  
Since $(p_{SO(2)})_* \co \pi_2 (K(\zed,2))\to \pi_2(BSO(2))$ is multiplication 
by 2, we have $2d(s,s')= c_1(s') - c_1(s)$. 
Equivalently, $c_1(s+a) = c_1(s)+2a$ for $a =d(s,s')$. 
The assertion about conjugation is clear from the way it lifts to 
$\Bspin^c(H,w)$.\qed

Now suppose we are given pairs $(BH,w)$ and $(BH',w')$ as before, and 
a map $h\co BH\to BH'$ covered by a bundle map 
$\tilde h\co \xi_H \to \xi_{H'}$, with $h^* w' = w$. 
Then any $H$--bundle $\xi\to X$ determines an $H'$--bundle $\xi'\to X$ 
with the same $w_2$ and $W_3$, and $h$ determines maps $\Bspin (H,w)\to 
\Bspin (H',w')$ and $\Bspin^c (H,w) \to \Bspin^c (H',w')$. 
We obtain canonical equivariant identifications $\cS(\xi,w)\cong \cS(\xi',w')$ 
and $\cS^c (\xi,w)\cong \cS^c(\xi',w')$, and the latter preserves Chern classes 
and conjugation. 
On the other hand, given an $H$--bundle map $\tilde g\co \xi_1\to \xi_2$ 
covering $g\co X_1\to X_2$, we have induced maps $g^* \co\cS(\xi_2,w)\to \cS(\xi_1,w)$ 
and $g^*\co \cS^c(\xi_2,w) \to \cS^c (\xi_1,w)$ that are equivariantly equivalent 
to $g^* \co H^1(X_2;\zed_2) \to H^1 (X_1;\zed_2)$ and $g^*\co H^2(X_2;\zed) 
\to H^2( X_1;\zed)$ when the domains are nonempty, and characteristic 
classes and conjugation are preserved in the obvious way. 
If $g$ is a homotopy equivalence, then the maps $g^*$ are isomorphisms. 

\rk{Examples 2}%
(a)\stdspace If $h\co BSO(m)\to BSO(n)$, $2\le m\le n$, is induced by the usual 
inclusion of 
groups, we recover the stabilization-invariance of classical spin and 
spin$^c$--structures. 
We are free to pass to the limiting group $SO$, eliminating the dependence 
on $n$.

(b)\stdspace An oriented topological $n$--manifold $X$ has the homotopy type of 
a $CW$--complex, and it has a tangent bundle classified by a map into the 
universal bundle over $BSTOP (n)$ (eg \cite{\KS}). 
There is a canonical map $h\co BSO(n)\to BSTOP(n)$ that corresponds to
interpreting $\xi_{SO(n)}$
as a topological bundle and is a $\pi_2$--isomorphism 
of simply connected spaces. 
We immediately obtain a theory of spin and spin$^c$--structures on oriented
topological manifolds by using their tangent bundles (stabilized if $n<2$). 
As before, the theory is stabilization-invariant, and we can pass to the 
limiting case of $BSTOP$. On smooth manifolds,
the new theory canonically reduces via $h$ to the classical theory.
However, any orientation-preserving homeomorphism $g\co X_1\to X_2$ induces
an isomorphism of topological tangent bundles, hence, isomorphisms 
$g^* \co\cS(X_2) \cong \cS(X_1)$ and $g^* \co\cS^c (X_2) \cong \cS^c (X_1)$ 
as above. 
\endrk

To generalize to homotopy equivalences, we need one further construction. 
Suppose we are given a bundle map 
$$\matrix{ 
\strt\xi_H \times \xi_{H'}& \buildrel {\tilde k}\over\longrightarrow 
&\xi_{H''}\cr 
\Big\downarrow&&\Big\downarrow\cr
BH\times BH'&\buildrel k\over\longrightarrow&BH''\cr}
$$ 
with $k^* (w'') = w+w'$. 
Then a pair of bundles $\xi,\xi'\to X$ classified by $BH,BH'$ 
determine an $H''$--bundle $\xi''\to X$, and $w_2$ and $W_3$ add.

\proclaim{Proposition 3}
{A trivialization of $\xi''$ induces equivariant isomorphisms
$k_*\co$\break$\cS(\xi,w) 
\to \cS(\xi',w')$ and $k_*\co\cS^c (\xi,w)\to \cS^c (\xi',w')$, and the latter 
preserves conjugation and Chern classes.} 
\endproc

\prf
By obstruction theory, the map $k$ uniquely determines
a map $\hat k$ making the diagram 
$$\matrix{
\strt\Bspin (H,w)\times\Bspin (H',w')&\buildrel {\hat k}\over\longrightarrow 
&\Bspin (H'',w'')\cr 
\bigg\downarrow {\scriptstyle p_1}&&\bigg\downarrow {\scriptstyle p_2}\cr
BH\times BH'&\buildrel k\over\longrightarrow &BH''\cr}
$$ 
commute, and a similar diagram is induced for {\rm spin}$^c$ via the map 
$k\times k_0$, where $k_0 \co BSO (2) \times BSO(2)\to BSO(2)$ induces addition 
on $\pi_2$. 
The diagrams determine a map $k_{\#} \co\cS(\xi,w)\times \cS(\xi',w')\to 
\cS(\xi'',w'')$ and similarly for $\cS^c$. 
In the latter case, $k_{\#}$ commutes with conjugation and adds Chern classes. 
In either case, $\hat k$ restricts to addition
on the homotopy groups of the fibers of 
$p_1$ and $p_2$, so difference classes add under $k_{\#}$, and for suitably 
chosen base points $k_{\#}$ is given by addition on $H^1(X;\zed_2)$ or 
$H^2(X;\zed)$ whenever its domain is nonempty. 
Now a trivialization of $\xi''$ determines a trivial spin$^c$--structure 
$s'' \in \cS^c(\xi'',w'')$. 
Since $W_3(\xi) + W_3(\xi') = W_3 (\xi'')=0$, it follows that 
$\cS^c(\xi,w)$ is nonempty if and only if $\cS^c (\xi',w')$ is. 
For each $s\in \cS^c(\xi,w)$ there is a unique ``inverse'' $s'\in \cS^c(\xi',w')$ 
with $k_{\#} (s,s')=s''$. 
Let $k_* (s)$ equal the conjugate of $s'$. 
Then $k_* \co\cS^c (\xi,w)\to \cS^c (\xi',w')$ is an equivariant isomorphism, 
and it preserves conjugation and Chern classes since $s''$ is 
conjugation-invariant with $c_1(s'')=0$. 
A similar procedure (with $k_*(s)=s'$) works for spin structures.\qed

\rk{Example 4}
Any oriented, smooth $n$--manifold $X$ admits a unique isotopy class of 
proper embeddings in $\real^N$ for $N$ sufficiently large. 
This determines a normal bundle $\nu X$ that is unique up to stabilization. 
Since the tangent bundle $\tau X$ satisfies $\tau X\oplus \nu X = \tau\real^N
|X$ and the latter bundle is canonically trivial, the obvious
map $BSO(n) \times 
BSO(N-n)\to BSO(N)$ determines canonical equivariant identifications 
$\cS(\nu X,w^*) \cong \cS(\tau X,w^*)$ and $\cS^c (\nu X,w^*)\cong \cS^c(\tau X,
w^*)$, the latter preserving Chern classes and conjugation. 

\proclaim{Theorem 5}
{Let $(X,\partial X)$ be an oriented, possibly noncompact Poincar\'e
pair.  There is a canonical procedure for defining sets $\cS(X)$ and
$\cS^c(X)$ of spin and spin$^c$--structures on $X$ having the
structure described in Proposition~\PropOne\ (with respect to the
usual classes $w_2(X)$ and $W_3(X)$).  For $(X,\partial X)$ a smooth
manifold, the theory is canonically equivariantly equivalent to the
standard one (preserving Chern classes and conjugation).  For pairs
$(X_i,\partial X_i)$ as above, any orientation-preserving, pairwise,
proper homotopy equivalence $g\co(X_1,\partial X_1)\to (X_2,\partial
X_2)$ induces equivariant isomorphisms $g^* \co\cS(X_2) \cong
\cS(X_1)$ and $g^* \co$\break$ \cS^c (X_2)\cong \cS^c (X_1)$, the latter
preserving Chern classes and conjugation, and the construction is
functorial for such maps $g$.}
\endproc

\prf
The pair $(X,\partial X)$ has a canonical {\it Spivak normal fibration\/} 
\cite{\S} defined by embedding $(X,\partial X)$ pairwise and properly 
in half-space $\real^N \times ([0,\infty),\{ 0\} )$
(uniquely for $N$ sufficiently large), and making 
a fibration out of the collapsing map of the boundary of a regular 
neighborhood. 
The resulting  oriented spherical fibration over $X$ is classified by a 
fiber-preserving map into the universal spherical fibration, whose base space 
stabilizes to $BSG$. 
As in Example~2(b), there is a canonical map $h\co BSO \to BSG$ induced by the
spherical fibrations $\xi_{SO(n)}-(0\hbox{--section})$, and $h$ is a 
$\pi_2$--isomorphism of simply connected spaces. 
We immediately obtain $\cS(X)$, $\cS^c(X)$ and characteristic classes 
satisfying Proposition~\PropOne, using the Spivak fibration and $BSG$. 
(The resulting classes $w_2(X)$ and $W_3(X)$ are well known.)
For $(X,\partial X)$ a smooth manifold, the theory is canonically equivalent 
(via $h$) to that of the stable normal bundle, which is the usual theory 
over the tangent bundle by Example~4. 
A homotopy equivalence $g$ as above induces a fiber-preserving map of the 
corresponding Spivak fibrations, and hence, the required maps $g^*$.\qed 
\ppar

The map $g^* \co\cS^c(X_2) \cong \cS^c(X_1)$ induced by a homotopy equivalence 
can frequently be computed explicitly. 
We consider the case where $X_2$ contains a 1--dimensional subcomplex 
with a regular neighborhood $N_2$ that is a manifold, such that $H^2(X_2,N_2;
\zed)$ has no 2--torsion. 
We also assume that $g\co X_1\to X_2$ restricts to a homeomorphism from 
$N_1 = g^{-1}(N_2)$ to $N_2$. 
These conditions are always satisfied if $g$ is a homeomorphism between smooth 
manifolds, for example by taking $N_2$ to be 
a neighborhood of the 1--skeleton of $X_2$. 
Now the map $g^* \co H^* (X_2,N_2) \cong H^* (X_1,N_1)$ is an isomorphism. 
A (stable) trivialization $\tau_2$ of the tangent bundle 
of $N_2$ (or equivalently, 
of the stable normal bundle) pulls back via $g|N_1$ to a trivialization 
$\tau_1$ over $N_1$, and $g^* w_2 (X_2,\tau_2) = w_2 (X_1,\tau_1)$. 
Given spin$^c$--structures $s_i \in \cS^c (X_i)$, pick any trivializations 
$\hat \tau_i$ of $s_i|N_i$ over $\tau_i$. 
Then by Proposition~\PropOne, 
$g^* c_1 (s_2,\hat\tau_2) - c_1(s_1,\hat\tau_1)$ reduces to zero mod~2.  
Since  $H^2(X_1,N_1;\zed)$ has no 2--torsion, there is a unique class 
$\delta (s_1,s_2)\in H^2 (X_1,N_1;\zed)$ with $2\delta (s_1,s_2) = 
g^* c_1 (s_2,\hat\tau_2) - c_1(s_1,\hat\tau_1)$. 
If we change $\hat\tau_i$ with $\tau_i$ fixed, then $\delta (s_1,s_2)$ changes 
by the coboundary of a cochain in $N_1$, so it represents a class $d(s_1,s_2)
\in H^2 (X_1;\zed)$ that depends only on $s_1$ and $s_2$ ($\tau_i$ fixed). 
But $\delta(s_1,s_2)$ vanishes for $s_1 = g^*s_2$ and $\hat\tau_1$
given by pulling back $\hat\tau_2$, and a change of $s_i$ changes
$2 \delta(s_1,s_2)$ by twice the corresponding relative difference class
(by the addition formula of Proposition~\PropOne\ applied to $X_i/N_i$).
Thus, $d(s_1,s_2)$ is precisely the difference class 
$d(s_1,g^* s_2)$, in a form accessible to computation. 

\rk{Remarks}%
(a)\stdspace Spin$^c$--structures have several other convenient characterizations.  
As we observed in proving Proposition~\PropOne, 
a spin$^c$--structure on $\xi\to X$ 
is the same as a line bundle $L$ and spin structure on $\xi\times L\to X$. 
For a different approach, recall that Milnor \cite{\Mil} observed that a spin 
structure on an oriented vector bundle over a $CW$--complex is equivalent (after 
stabilizing if necessary) to a trivialization over the 1--skeleton that can 
be extended over the 2--skeleton, just as an orientation is a trivialization 
over the 0--skeleton that extends over the 1--skeleton. 
Similarly, a spin$^c$--structure over an oriented vector bundle is 
equivalent (after stabilizing if the fiber dimension is odd or $\le2$)
to a complex 
structure over the 2--skeleton that can be extended over the 3--skeleton. 
To see this, observe that the map of classifying spaces induced by
inclusion $i\co U(n)\to SO(2n)$ lifts canonically to a map $j\co BU(n) \to 
B\hbox{spin}^c (SO(2n),w^*)$ by first lifting the map $id\times B\det :BU(n) 
\to BU(n)\times BSO(2)$ to $\Bspin^c (U(n),i^* w^*)$.  (In fact, the 
corresponding diagram exists on the group level.)
Thus, any complex structure determines a spin$^c$--structure 
(and the correspondence preserves $c_1$ and conjugation). For $n \ge 2$, 
this correspondence is bijective for 2--complexes and surjective for 
3-complexes, since the map $j$ has a 2--connected fiber.
The observation now follows from the fact that restriction induces a 
bijection from spin$^c$--structures to those over the 2--skeleton extending 
over the 3--skeleton.  The same remark applies to bundles classified by
$BSTOP$ or $BSG$ if we define a complex structure to be a lift of the
classifying map to $BU$.

(b)\stdspace The Wu relations are known to hold for Poincar\'e complexes. 
In particular, for a compact, oriented 4--dimensional Poincar\'e complex $X$ 
(without boundary) we have $w_2(X)\cup x = x\cup x$ for all $x\in H^2(X;
\zed_2)$. 
The usual argument \cite{\HH} then shows that $W_3(X)=0$, so all such 
complexes admit spin$^c$--structures. 

(c)\stdspace As in the classical case, we have a canonical map $i:\Bspin (H,w)\to 
\Bspin^c (H,w)$ as the fiber of $p_{SO(2)}$ (induced by inclusion of groups), 
inducing a map $\alpha\,\co\cS(\xi,w)\to \cS^c (\xi,w)$ that is equivariantly  
equivalent (when the domain is nonempty) 
to the Bockstein homomorphism $\beta \co H^1(X;\zed_2) \to H^2(X;\zed)$. 
The image $\hbox{Im } \alpha$ is the set of
spin$^c$--structures with $c_1=0$, or equivalently, 
the set of conjugation-invariant structures.
To verify that $\alpha$ has the stated equivariance and image, note that
we can either consider $i$ to be an inclusion into the fixed set of
conjugation or replace it by a fibration $p$. Over each point in $BH$,
$i$ and $p$ will restrict to the canonical inclusion and fibration
$\real P^\infty \to \cee P^\infty$, respectively, both of which represent
the unique nontrivial homotopy class of maps in
$[\real P^\infty, \cee P^\infty]$. For a fixed classifying map
$\tilde f\co\xi \to \xi_H$, spin structures $s_1,s_2 \in \cS(\xi,w)$ determine
lifts $\hat f_1,\hat f_2\co X \to \Bspin (H,w)$. We can assume that these
agree over the 0--skeleton and that $p \circ \hat f_1$, $p\circ\hat f_2$
agree over the 1--skeleton, giving us obstruction cochains $d(s_1,s_2) \in
C^1(X;\zed_2)$ and $d(\alpha s_1,\alpha s_2)\in C^2(X;\zed)$.  Now 
$d(\alpha s_1, \alpha s_2)$ evaluated on a 2--cell $c$ is the element of
$\pi_2(\cee P^\infty) \cong \zed$ given by $p \circ \hat f_2(c) -
p\circ \hat f_1(c)$.  Since the boundary operator $\pi_2(\cee P^\infty)
\to \pi_1(\real P^1)$ of $p$ is multiplication by 2, the same coefficient
is obtained as $\frac12\langle d(s_1,s_2),\partial c\rangle
= \langle\beta d(s_1,s_2),c\rangle$.
Thus, we obtain the required equivariance $d(\alpha s_1, \alpha s_2) =
\beta d(s_1,s_2)$.  To compute $\hbox{Im }\alpha$, first note that any
$s \in \hbox{Im }\alpha$ is conjugation-invariant (since $i$ is) with $c_1=0$.
If $\cS (\xi,w)$ is nonempty, fix $s \in\hbox{Im }\alpha$ and let $s'$ be any
spin$^c$--structure that either is conjugation-invariant or satisfies 
$c_1(s')=0$. By Proposition~\PropOne, $2d(s,s')=0$, so $d(s,s') \in\hbox{Im }\beta$
and $s' \in\hbox{Im }\alpha$.  It now suffices to show that when $\cS (\xi,w)$
is empty, no spin$^c$--structure has $c_1=0$ or is conjugation-invariant.
The first assertion is obvious since $c_1|_2=w_2 \ne 0$.  For the remaining
assertion, choose $s \in \cS^c(\xi,w)$ with conjugate $\bar s$.  Since
$\pi_1(\cee P^\infty,\real P^\infty) =0$, we can assume that the lift
$\hat f\co X \to \Bspin^c(H,w)$ determined by $s$ maps the 1--skeleton $X_1$
into $i(\Bspin (H,w))$, which is fixed by conjugation. Thus, $\hat f$ and
its conjugate determine a difference cochain $d(s,\bar s) \in C^2(X,\zed)$.
Since $\pi_2(\cee P^\infty)\to
\pi_2(\cee P^\infty,\real P^\infty)$ is multiplication by 2 on $\zed$,
we can change $d(s,\bar s)$ by any coboundary by changing 
$\hat f|X_1\co X_1 \to i(\Bspin (H,w))$.  Thus, if $s=\bar s$ we can assume
that $d(s,\bar s)=0$, so over each 2--cell, $\hat f$ is conjugation-invariant
up to homotopy rel $\partial$.  But conjugation fixes only 0 in
$\pi_2(\cee P^\infty,\real P^\infty)$, so $\hat f$ can then be homotoped
into $i(\Bspin (H,w))$, ie $s \in \hbox{Im }\alpha$.
%
\vglue0.2truein
\references         

\end